\documentclass[12pt,reqno]{amsart}
\usepackage{amsfonts}
\usepackage{amssymb, amsmath, latexsym}
\usepackage{lipsum}
\usepackage{hyperref}
\usepackage{longtable}
\usepackage[square,numbers]{natbib}
\usepackage{courier}

\usepackage{mathtools}
\DeclarePairedDelimiter\floor{\lfloor}{\rfloor}
\usepackage{array}

\numberwithin{equation}{section}
\theoremstyle{plain}
\newtheorem{theorem}{Theorem}[section]

\newcommand\T{\rule{0pt}{2.6ex}}
\newcommand\B{\rule[-1.2ex]{0pt}{0pt}}

\title[Witness identities for three Ramanujan congruences] {Witness identities for three Ramanujan congruences}
\author[N. D. Baruah]{Nayandeep Deka Baruah}
\address{Department of Mathematical Sciences, Tezpur University, Napaam 784028, Sonitpur, Assam, India}
\email{nayan@tezu.ernet.in}
\author[H. Das]{Hirakjyoti Das}
\address{Department of Mathematical Sciences, Tezpur University, Napaam 784028, Sonitpur, Assam, India}
\email{hdas@tezu.ernet.in}
\author[A. Sarma]{Abhishek Sarma}
\address{Department of Mathematical Sciences, Tezpur University, Napaam 784028, Sonitpur, Assam, India}
\email{msp21107@tezu.ac.in}

\allowdisplaybreaks
\begin{document}
\begin{abstract}
For the unrestricted partition function $p(n)$ for integers $n\ge0$, it is known that $p(49n+19)\equiv 0 \pmod{49}$, $p(49n+33)\equiv 0 \pmod{49}$, and $p(49n+40)\equiv 0 \pmod{49}$ for all $n\ge0$. We find witness identities for these Ramanujan congruences.
\end{abstract}
\maketitle
\noindent{\footnotesize Key words:   Partition function, Generating function, Congruence}

\vskip 3mm
\noindent {\footnotesize 2010 Mathematical Reviews Classification Numbers: 05A17, 11P83}

\section{Introduction}
For complex numbers, $a$ and $q$ with $\mid q\mid<1$, we define the standard infinite $q$-product as 
\begin{align*}
 \left(a;q\right)_\infty:=\prod_{j=0}^{\infty}\left(1-a q^j\right).
\end{align*}
Throughout the paper, for all $j\ge 1$, $k\ge 1$, and $n\ge1$, we adopt 
\begin{align}
f_n&:=\left(q^n;q^n\right)_\infty,\notag\\
\label{f_j,k} f_{j,k}(q)&:=\dfrac{\left(q^{2j};q^k\right)_\infty\left(q^{k-2j};q^k\right)_\infty}{\left(q^{j};q^k\right)_\infty\left(q^{k-j};q^k\right)_\infty}.
\end{align}

A partition $\lambda:=(\lambda_1, \lambda_2,\ldots,\lambda_k)$ of a positive integer $n$ is a finite non-increasing sequence of positive integers $\lambda_1, \lambda_2,\ldots,\lambda_k$ such that $\lambda_1+\lambda_2+\cdots+\lambda_k=n$. We call the $\lambda_j$'s as the parts of $\lambda$. The study of partition theory has a long history since Euler (1707--1783) gave the generating function for the number of partitions of a positive integer $n$. If $p(n)$ counts the partitions of $n$, then its generating function is given by
\begin{align*}
    \sum_{n=0}^\infty p(n)q^n=\dfrac{1}{f_1},
\end{align*}
where, by convention, $p(0)=1$. The arithmetic properties of the partition function $p(n)$ have been studied extensively  after Ramanujan \cite{ram1, ram2, ram3} found his famous congruences modulo 5, 7, and 11, namely, for all $n\ge 0$,
\begin{align}
\label{Rama Cong 1} p(5n+4)&\equiv 0~(\textup{mod}~5),\\
\label{Rama Cong 2}p(7n+5)&\equiv 0~(\textup{mod}~7),\\
\label{Rama Cong 3}p(11n+6)&\equiv 0~(\textup{mod}~11).
\end{align}
Ramanujan \cite{ram1} also conjectured the following infinite family of congruences. If $\delta=5^a7^b11^c$ and $\gamma$ is an integer such that $24\gamma\equiv 1~(\textup{mod}~\delta)$, then for all $n\ge0$,
\begin{align}
\label{Ram Conj 1} p\left(\delta n+\gamma\right) \equiv 0~~ \left(\textup{mod}~\delta\right).
\end{align}
In \cite{BK99}, Ramanujan proved \eqref{Ram Conj 1}  when $a$ is arbitrary and $b=c=0$. Atkin \cite{Atkin67} provided a proof of \eqref{Ram Conj 1} for  arbitrary $c$ and $a=b=0$. However   when $b$ is arbitrary and $a=c=0$, \eqref{Ram Conj 1} turned out to be incorrect. The corrected version of \eqref{Ram Conj 1} is stated below. 

If $\delta^{\prime} =5^a 7^{b^\prime} 11^c$, where $b^\prime=b$ when $b=0,1,2,$ and $b^\prime=\floor*{(b+2)/2}$ when $b>2$, and $\gamma$ is an integer such that $24\gamma\equiv 1~(\textup{mod}~\delta)$, then for all $n\ge0$,
\begin{align}
 \label{Ram Conj 1 Corrected} p\left(\delta n+\gamma\right) \equiv 0~~ \left(\textup{mod}~\delta^\prime\right).
\end{align}
The case when $b$ is arbitrary and $a=c=0$ in the above corrected version was proved by Watson \cite{Watson38}. 

Ramanujan also discovered the following generating functions, from which \eqref{Rama Cong 1} and \eqref{Rama Cong 2} follow immediately.
\begin{align}
\sum_{n=0}^\infty p(5n+4)q^n &= 5\dfrac{f_5^5}{f_1^6}, \notag\\
\label{Witness 7n+5} \sum_{n=0}^\infty p(7n+5)q^n &= 7\dfrac{f_7^3}{f_1^4} + 49 q\dfrac{f_7^7}{f_1^8}.
\end{align}
Such generating functions are called \textit{witness identities} for Ramanujan congruences. A witness identity for \eqref{Rama Cong 3} was first found by Lehner \cite{Lehner43} and later by Atkin \cite{Atkin67}, Paule and Radu \cite{PR17}, and Goswami, Jha, and Singh \cite{GJS22}. Hirschhorn and Hunt \cite{HH1} later found a witness identity for the infinite family \eqref{Ram Conj 1} for arbitrary $a$ and $b=c=0$. A witness identity for the family \eqref{Ram Conj 1 Corrected} when $b$ is arbitrary and $a=c=0$ was found by Garvan \cite{Garvan84}.

In \cite{ram2}, Ramanujan mentioned that for all $n\ge0$,
\begin{align}
\label{Mod 49 19} p(49n+19)&\equiv 0 \pmod{49},\\
\label{Mod 49 33}p(49n+33)&\equiv 0 \pmod{49},\\
\label{Mod 49 40}p(49n+40)&\equiv 0 \pmod{49},
\end{align}
which do not come under \eqref{Ram Conj 1 Corrected}. In 1980, an infinite family of congruences modulo powers of 7 was found by Ramanathan \cite[Corollary 1]{Ramanathan80}, which generalizes \eqref{Mod 49 19}--\eqref{Mod 49 40}.

A simple proof of \eqref{Mod 49 19}--\eqref{Mod 49 40} is as follows. 
One of Jacobi's identities is \cite[Theorem 1.3.9]{Spirit}
\begin{align}
\label{Jacobi} f_1^3 &=\sum_{n=0}^\infty (-1)^n(2n+1)q^{n(n+1)/2}.
\end{align} 
Taking modulo 49 in \eqref{Witness 7n+5} and then using \eqref{Jacobi}, we have
\begin{align}
\label{Step 1}\sum_{n=0}^\infty p(7n+5)q^n &= 7\dfrac{f_7^3}{f_1^4} + 49 q\dfrac{f_7^7}{f_1^8}\equiv 7\dfrac{f_7^3}{f_1^4} \equiv 7f_1^3f_7^2 \notag\\
&\equiv 7f_7^2\sum_{n=0}^\infty (-1)^n(2n+1)q^{n(n+1)/2}\pmod{49}.
\end{align} 
It can be easily seen that $n(n+1)/2\not\equiv 2,4,5\pmod{7}$ for all $n\ge0$. Therefore, from \eqref{Step 1}, extracting the terms involving $q^{7n+2}$, $q^{7n+4}$, and $q^{7n+5}$ for all $n\ge0$, we obtain \eqref{Mod 49 19}--\eqref{Mod 49 40}.

In this paper, we find three witness identities for \eqref{Mod 49 19}--\eqref{Mod 49 40}, which we state in the following theorem. Note that the huffing operator method in \cite{Garvan84} by Garvan does not provide the arithmetic progressions of  \eqref{Mod 49 19}--\eqref{Mod 49 40}.
\begin{theorem}\label{Theorem mod 49} 
We have
\begin{align}
\label{Witness 49n+19} \sum_{n=0}^\infty p\left(49n+19\right)q^n=
~&49\dfrac{f_{7}^{28}}{f_{1}^{29}}\Bigg(\sum_{j=0}^{16}\alpha_{19,j} f_{1,7}^{23-j}(q)f_{2,7}^{2j}(q)+\sum_{j=0}^{2}\beta_{19,j} f_{1,7}^{3j+26}(q)f_{3,7}^{2j+2}(q)\notag\\
&+\dfrac{f_{7}^{28}}{f_{1}^{28}}\Bigg(\sum_{j=0}^{33}\gamma_{19,j} f_{1,7}^{47-j}(q)f_{2,7}^{2j+1}(q)+\sum_{j=0}^{6}\delta_{19,j} f_{1,7}^{3j+49}(q)f_{3,7}^{2j+1}(q)
\Bigg)\Bigg),\\
\label{Witness 49n+33}\sum_{n=0}^\infty p\left(49n+33\right)q^n=
~&49\dfrac{f_{7}^{28}}{f_{1}^{29}}\Bigg(\sum_{j=0}^{16}\alpha_{33,j} f_{1,7}^{22-j}(q)f_{2,7}^{2j}(q)+\sum_{j=0}^{2}\beta_{33,j} f_{1,7}^{3j+25}(q)f_{3,7}^{2j+2}(q)\notag\\
&+\dfrac{f_{7}^{28}}{f_{1}^{28}}\Bigg(\sum_{j=0}^{33}\gamma_{33,j} f_{1,7}^{46-j}(q)f_{2,7}^{2j+1}(q)+\sum_{j=0}^{6}\delta_{33,j} f_{1,7}^{3j+48}(q)f_{3,7}^{2j+1}(q)
\Bigg)\Bigg),\\
\label{Witness 49n+40}\sum_{n=0}^\infty p\left(49n+40\right)q^n=
~&49\dfrac{f_{7}^{28}}{f_{1}^{29}}\Bigg(\sum_{j=0}^{15}\alpha_{40,j} f_{1,7}^{21-j}(q)f_{2,7}^{2j+1}(q)+\sum_{j=0}^{3}\beta_{40,j} f_{1,7}^{3j+23}(q)f_{3,7}^{2j+1}(q)\notag\\
&+\dfrac{f_{7}^{28}}{f_{1}^{28}}\Bigg(\sum_{j=0}^{33}\gamma_{40,j} f_{1,7}^{46-j}(q)f_{2,7}^{2j}(q)+\sum_{j=0}^{6}\delta_{40,j} f_{1,7}^{3j+49}(q)f_{3,7}^{2j+2}(q)
\Bigg)\Bigg),
\end{align}
where $\alpha_{r,j}$,  $\beta_{r,j}$, $\gamma_{r,j}$, and $\delta_{r,j}$ for $r\in\{19,33,40\}$ are given in Appendix.
\end{theorem}


\section{Proof of Theorem \ref{Theorem mod 49}}
The proofs of \eqref{Witness 49n+19}--\eqref{Witness 49n+40} are similar. Therefore, we prove \eqref{Witness 49n+19} only in detail. We let
\begin{align*}
\sum_{n=0}^\infty p_1(n)q^n := \dfrac{f_7^3}{f_1^4} \quad\text{and}\quad \sum_{n=0}^\infty p_2(n)q^n := q\dfrac{f_7^7}{f_1^8}
\end{align*}
so that by \eqref{Witness 7n+5}, we have
\begin{align}
\label{Witness 49n+19 Eve 1} \sum_{n=0}^\infty p(7n+5)q^n &=7\sum_{n=0}^\infty p_1(n)q^n+49\sum_{n=0}^\infty p_2(n)q^n.
\end{align}
Now, to prove \eqref{Witness 49n+19}, we require the exact generating functions for $p_1(7n+2)$ and $p_2(7n+2)$ for all $n\ge0$, both of which can be found in a similar way. Therefore, we find the exact generating function for  $p_1(7n+2)$ only.

Firstly, we require the following 7-dissection of $f_1$ \cite[p. 274, Theorem 12.1 for $n=7$]{B1},
\begin{align}
\label{7 Dissection f1}f_1&=f_{49}\left(f_{1,7}\left(q^7\right)-q f_{2,7}\left(q^7\right)-q^2+q^5f_{3,7}\left(q^7\right)\right),
\end{align}
where $f_{j,k}(q)$ are defined in \eqref{f_j,k}. Using Jacobi's identity \eqref{Jacobi} and \eqref{7 Dissection f1}, we have
\begin{align*}
\sum_{n=0}^\infty (-1)^n(2n+1)q^{n(n+1)/2}&=\left( f_{49}\left(f_{1,7}\left(q^7\right)-q f_{2,7}\left(q^7\right)-q^2+q^5f_{3,7}\left(q^7\right)\right)\right)^3.
\end{align*}
Expanding the right side and using the fact that $n(n+1)/2\not\equiv 2,4,5\pmod{7}$ for all $n\ge0$, we easily arrive at the following equalities.
\begin{align}
\label{Condition 1} f_{1,7}^2\left(q^7\right)-f_{1,7}\left(q^7\right) f_{2,7}^2\left(q^7\right)-q^7f_{3,7}\left(q^7\right)&=0,\\
f_{1,7}\left(q^7\right)-f_{2,7}^2\left(q^7\right)-q^7 f_{2,7}\left(q^7\right) f_{3,7}^2\left(q^7\right)&=0,\\
f_{2,7}\left(q^7\right)-f_{1,7}^2\left(q^7\right) f_{3,7}\left(q^7\right)+q^7 f_{3,7}^2\left(q^7\right)&=0,\\
\label{Condition 2} f_{1,7}\left(q^7\right) f_{2,7}\left(q^7\right) f_{3,7}\left(q^7\right)&=1.
\end{align}

Secondly, we need the following binomial theorem. For positive integers $\ell\ge1$, we have
\begin{align}
\label{Binomial Theorem} (a+b+c+d)^\ell=\sum_{i=0}^{\ell}\sum_{j=0}^{\ell-i}\sum_{k=0}^{i}\binom{\ell}{i}\binom{\ell-i}{j}\binom{i}{k}a^{\ell-i-j}b^{j}c^{i-k}d^{k}.
\end{align}

Now,
\begin{align}
\label{Witness 49n+19 Eve 1.5} \sum_{n=0}^\infty p_1(n)q^n&=\dfrac{f_7^3}{f_1^4}=\left(q^7;q^7\right)_\infty^3\dfrac{\displaystyle{\prod_{j=1}^{6}}\left(\zeta^j q;\zeta^j q\right)_\infty^4}{\displaystyle{\prod_{j=0}^{6}}\left(\zeta^j q;\zeta^j q\right)_\infty^4},
\end{align}
where $\zeta\neq1$ is a seventh root of unity. Let $\displaystyle{\prod_{j=1}^{6}}\left(\zeta^j q;\zeta^j q\right)_\infty^4=:N(q)$ and $\displaystyle{\prod_{j=0}^{6}}\left(\zeta^j q;\zeta^j q\right)_\infty=:D(q)$. Then by the above identity,
\begin{align}
\label{Witness 49n+19 Eve 2} D^4(q)=f_1^4\times N(q)
\end{align}
and by the fact that $\displaystyle{\prod_{j=0}^6\left(1-\zeta^jq\right)=\left(1-q^7\right)}$, we have
\begin{align}
\label{Witness 49n+19 Eve 3}D(q)=\dfrac{f_7^8}{f_{49}}.
\end{align}

Using \eqref{7 Dissection f1}, \eqref{Binomial Theorem} with $a=f_{1,7}\left(q^7\right)$, $b=-qf_{2,7}\left(q^7\right)$, $c=-q^2$, $d=q^5f_{3,7}\left(q^7\right)$, and $\ell=4$, and $\displaystyle{N(q)=\sum_{j=0}^{6}N_{j}(q)=:\sum_{j=0}^{6}\sum_{n=0}^{\infty}N_j(7n+j)q^{7n+j}}$ in \eqref{Witness 49n+19 Eve 2}, we find that
\begin{align}
\label{Expansion 1} \dfrac{D^4(q)}{f_{49}^4}&=\Bigg( \sum_{i=0}^{4}\sum_{j=0}^{4-i}\sum_{k=0}^{i}\binom{4}{i}\binom{4-i}{j}\binom{i}{k}(-1)^{i+j-k}\notag\\
&\quad\quad f_{1,7}^{4-i-j}\left(q^7\right)f_{2,7}^j\left(q^7\right)f_{3,7}^k\left(q^7\right)q^{2i+j+3k}\Bigg)\times \left( \sum_{\ell=0}^{6}N_{\ell}(q) \right).
\end{align}
Now, after expanding the right side of the above identity, we rearrange it according to the residues of the exponents of $q$ modulo 7 as follows. For that, we first explicitly mention the  powers $q^{2i+j+3k}$ for different values of $i$, $j$, and $k$ in the following table. 

\begin{longtable}{c c c c c c c}
\caption{Explicit tuples $\left(i,j,k,q^{2i+j+3k}\right)$}
\label{Explicit i j k} 
\\
\hline
$\left(0,0,0,q^0\right)$ & $\left(0,1,0,q\right)$  & $\left(0,2,0,q^2\right)$ & $\left(0,3,0,q^3\right)$ & $\left(0,4,0,q^4\right)$ & $\left(1,3,0,q^5\right)$ & $\left(1,1,1,q^6\right)$ \T\B\\
$\left(1,2,1,q^7\right)$ & $\left(1,3,1,q^8\right)$  & $\left(1,0,0,q^2\right)$ & $\left(1,1,0,q^3\right)$ & $\left(1,2,0,q^4\right)$ & $\left(1,0,1,q^5\right)$ & $\left(2,2,0,q^6\right)$ \T\B\\
$\left(2,0,1,q^7\right)$  & $\left(2,1,1,q^8\right)$  & $\left(2,2,1,q^9\right)$ & $\left(2,0,2,q^{10}\right)$ & $\left(2,0,0,q^4\right)$ & $\left(2,1,0,q^5\right)$ & $\left(3,0,0,q^6\right)$ \T\B\\
$\left(3,1,0,q^7\right)$    &  $\left(4,0,0,q^8\right)$  & $\left(3,0,1,q^9\right)$ & $\left(3,1,1,q^{10}\right)$ & $\left(2,1,2,q^{11}\right)$ & $\left(2,2,2,q^{12}\right)$ & $\left(3,1,2,q^{13}\right)$ \T\B\\
$\left(4,0,2,q^{14}\right)$      &  $\left(3,0,3,q^{15}\right)$ & $\left(3,1,3,q^{16}\right)$ & $\left(4,0,3,q^{17}\right)$ & $\left(4,0,1,q^{11}\right)$ & $\left(3,0,2,q^{12}\right)$ & $\left(4,0,4,q^{20}\right)$ \T\B\\
\hline
\end{longtable}
We set
\begin{align}
\label{R i j k} R_{i,j,k}(q):=\binom{4}{i}\binom{4-i}{j}\binom{i}{k}(-1)^{i+j-k} f_{1,7}^{4-i-j}\left(q^7\right)f_{2,7}^j\left(q^7\right)f_{3,7}^k\left(q^7\right)q^{2i+j+3k}.
\end{align}

Using Table \ref{Explicit i j k} and $R_{i,j,k}(q)$, \eqref{Expansion 1} can be rearranged as
\begin{align*}
\dfrac{D^4(q)}{E_{49}^4}&=\left( \sum_{i=0}^{4}\sum_{j=0}^{4-i}\sum_{k=0}^{i}R_{i,j,k}(q)\right)\times \left( \sum_{\ell=0}^{6}N_{\ell}(q) \right)\\
&=\Bigg(\Big(R_{0,0,0}(q)+R_{1,2,1}(q)+R_{2,0,1}(q)+R_{3,1,0}(q)+R_{4,0,2}(q)\Big) \\
&\quad +\Big(R_{0,1,0}(q)+R_{1,3,1}(q)+R_{2,1,1}(q)+R_{4,0,0}(q)+R_{3,0,3}(q)\Big)\\
&\quad +\Big(R_{0,2,0}(q)+R_{1,0,0}(q)+R_{2,2,1}(q)+R_{3,0,1}(q)+R_{3,1,3}(q)\Big)\\
&\quad +\Big(R_{0,3,0}(q)+R_{1,1,0}(q)+R_{2,0,2}(q)+R_{3,1,1}(q)+R_{4,0,3}(q)\Big)\\
&\quad +\Big(R_{0,4,0}(q)+R_{1,2,0}(q)+R_{2,0,0}(q)+R_{2,1,2}(q)+R_{4,0,1}(q)\Big)\\
&\quad +\Big(R_{1,3,0}(q)+R_{1,0,1}(q)+R_{2,1,0}(q)+R_{2,2,2}(q)+R_{3,0,2}(q)\Big)\\
&\quad +\Big(R_{1,1,1}(q)+R_{2,2,0}(q)+R_{3,0,0}(q)+R_{3,1,2}(q)+R_{4,0,4}(q)\Big)\Bigg)\times \left( \sum_{\ell=0}^{6}N_{\ell}(q) \right).
\end{align*}
Since $D(q)=f_7^8/f_{49}$, therefore the left side of the above equation is a function of $q^7$. So, the above equation can be put in the matrix form
\begin{align}
\label{Matrix Eq 1} M=A\times N,
\end{align}
where
\begin{align*}
M:=
\begin{pmatrix}
D^4(q)/f_{49}^4\\
0\\
0\\
0\\
0\\
0\\
0
\end{pmatrix},
\quad
N:=
\begin{pmatrix}
N_{0}(q)\\
N_{1}(q)\\
N_{2}(q)\\
N_{3}(q)\\
N_{4}(q)\\
N_{5}(q)\\
N_{6}(q)
\end{pmatrix},
\end{align*}
\begin{align*}
A:=
\begin{pmatrix}
A_{1}(q) & A_{7}(q) & A_{6}(q) & A_{5}(q) & A_{4}(q) & A_{3}(q) & A_{2}(q)\\
A_{2}(q) & A_{1}(q) & A_{7}(q) & A_{6}(q) & A_{5}(q) & A_{4}(q) & A_{3}(q)\\
A_{3}(q) & A_{2}(q) & A_{1}(q) & A_{7}(q) & A_{6}(q) & A_{5}(q) & A_{4}(q)\\
A_{4}(q) & A_{3}(q) & A_{2}(q) & A_{1}(q) & A_{7}(q) & A_{6}(q) & A_{5}(q)\\
A_{5}(q) & A_{4}(q) & A_{3}(q) & A_{2}(q) & A_{1}(q) & A_{7}(q) & A_{6}(q)\\
A_{6}(q) & A_{5}(q) & A_{4}(q) & A_{3}(q) & A_{2}(q) & A_{1}(q) & A_{7}(q)\\
A_{7}(q) & A_{6}(q) & A_{5}(q) & A_{4}(q) & A_{3}(q) & A_{2}(q) & A_{1}(q)
\end{pmatrix},
\end{align*}
where
\begin{align}
\label{A1 in f}A_1(q)&:=R_{0,0,0}(q)+R_{1,2,1}(q)+R_{2,0,1}(q)+R_{3,1,0}(q)+R_{4,0,2}(q),\\
\label{A2 in f}A_2(q)&:=R_{0,1,0}(q)+R_{1,3,1}(q)+R_{2,1,1}(q)+R_{4,0,0}(q)+R_{3,0,3}(q),\\
\label{A3 in f}A_3(q)&:=R_{0,2,0}(q)+R_{1,0,0}(q)+R_{2,2,1}(q)+R_{3,0,1}(q)+R_{3,1,3}(q),\\
\label{A4 in f}A_4(q)&:=R_{0,3,0}(q)+R_{1,1,0}(q)+R_{2,0,2}(q)+R_{3,1,1}(q)+R_{4,0,3}(q),\\
\label{A5 in f}A_5(q)&:=R_{0,4,0}(q)+R_{1,2,0}(q)+R_{2,0,0}(q)+R_{2,1,2}(q)+R_{4,0,1}(q),\\
\label{A6 in f}A_6(q)&:=R_{1,3,0}(q)+R_{1,0,1}(q)+R_{2,1,0}(q)+R_{2,2,2}(q)+R_{3,0,2}(q),\\
\label{A7 in f}A_7(q)&:=R_{1,1,1}(q)+R_{2,2,0}(q)+R_{3,0,0}(q)+R_{3,1,2}(q)+R_{4,0,4}(q).
\end{align}

Next, our aim is to find $N_{2}(q)$ \big(for proving \eqref{Witness 49n+33} and \eqref{Witness 49n+40}, we require to find $N_{4}(q)$ and $N_{5}(q)$, respectively\big) by solving \eqref{Matrix Eq 1}. Since all the entries of the column matrix $M$,  except the one in the first row, are zero, therefore, to find out  $N_{2}(q)$ \big($N_{4}(q)$ and $N_{5}(q)$\big), it enough to calculate the third  \big(the fifth and the sixth entries\big) entry in the first column of the inverse of the matrix $A$. For that, we first find the cofactor of the third \big(the fifth and the sixth entries\big) entry of the first row of the matrix $A$.
 
We accomplish the above task with the determinant function of a square matrix of order 5.  Let $M:=\left(m_{i,j}\right)_{1\le i, j\le 5}$ and 
\begin{align*}
\mathcal{A} &:=m_{3,3} \left(m_{4,4} m_{5,5}-m_{5,4} m_{4,5}\right)-m_{3,4} \left(m_{4,3} m_{5,5}-m_{5,3} m_{4,5}\right)+m_{3,5} \left(m_{4,3} m_{5,4}-m_{5,3} m_{4,4}\right),\\
\mathcal{B} &:=m_{3,2} \left(m_{4,4} m_{5,5}-m_{5,4} m_{4,5}\right)-m_{3,4} \left(m_{4,2} m_{5,5}-m_{5,2} m_{4,5}\right)+m_{3,5} \left(m_{4,2} m_{5,4}-m_{5,2} m_{4,4}\right),\\
\mathcal{C}&:=m_{3,2} \left(m_{4,3} m_{5,5}-m_{5,3} m_{4,5}\right)-m_{3,3} \left(m_{4,2} m_{5,5}-m_{5,2} m_{4,5}\right)+m_{3,5} \left(m_{4,2} m_{5,3}-m_{5,2} m_{4,3}\right),\\
\mathcal{D}&:=m_{3,2} \left(m_{4,3} m_{5,4}-m_{5,3} m_{4,4}\right)-m_{3,3} \left(m_{4,2} m_{5,4}-m_{5,2} m_{4,4}\right)+m_{3,4} \left(m_{4,2} m_{5,3}-m_{5,2} m_{4,3}\right),\\
\mathcal{E}&:=m_{3,1} \left(m_{4,4} m_{5,5}-m_{5,4} m_{4,5}\right)-m_{3,4} \left(m_{4,1} m_{5,5}-m_{5,1} m_{4,5}\right)+m_{3,5} \left(m_{4,1} m_{5,4}-m_{5,1} m_{4,4}\right),\\
\mathcal{F}&:= m_{3,1} \left(m_{4,3} m_{5,5}-m_{5,3} m_{4,5}\right)-m_{3,3} \left(m_{4,1} m_{5,5}-m_{5,1} m_{4,5}\right)+m_{3,5} \left(m_{4,1} m_{5,3}-m_{5,1} m_{4,3}\right),\\
\mathcal{G}&:= m_{3,1} \left(m_{4,3} m_{5,4}-m_{5,3} m_{4,4}\right)-m_{3,3} \left(m_{4,1} m_{5,4}-m_{5,1} m_{4,4}\right)+m_{3,4} \left(m_{4,1} m_{5,3}-m_{5,1} m_{4,3}\right),\\
\mathcal{H}&:= m_{3,1} \left(m_{4,2} m_{5,5}-m_{5,2} m_{4,5}\right)-m_{3,2} \left(m_{4,1} m_{5,5}-m_{5,1} m_{4,5}\right)+m_{3,5} \left(m_{4,1} m_{5,2}-m_{5,1} m_{4,2}\right),\\
\mathcal{I}&:= m_{3,1} \left(m_{4,2} m_{5,4}-m_{5,2} m_{4,4}\right)-m_{3,2} \left(m_{4,1} m_{5,4}-m_{5,1} m_{4,4}\right)+m_{3,4} \left(m_{4,1} m_{5,2}-m_{5,1} m_{4,2}\right),\\
\mathcal{J}&:= m_{3,1} \left(m_{4,2} m_{5,3}-m_{5,2} m_{4,3}\right)-m_{3,2} \left(m_{4,1} m_{5,3}-m_{5,1} m_{4,3}\right)+m_{3,3} \left(m_{4,1} m_{5,2}-m_{5,1} m_{4,2}\right).
\end{align*}
Then, it can be verified that
\begin{align}
\label{Det function}\textup{Det}(M)&=m_{1,1} \left(m_{2,2} \times \mathcal{A}- m_{2,3}\times \mathcal{B}+ m_{2,4} \times \mathcal{C}-  m_{2,5} \times \mathcal{D}\right)\notag\\
&\quad -m_{1,2} \left(  m_{2,1}\times \mathcal{A}-  m_{2,3}\times \mathcal{E}+  m_{2,4}\times \mathcal{F}-  m_{2,5}\times \mathcal{G}\right)\notag\\
&\quad +m_{1,3} \left(  m_{2,1}\times \mathcal{B} - m_{2,2} \times \mathcal{E}+  m_{2,4}\times \mathcal{H}-  m_{2,5} \times \mathcal{I}\right)\notag\\
&\quad -m_{1,4} \left(  m_{2,1} \times \mathcal{C}-  m_{2,2} \times \mathcal{F}+  m_{2,3} \times \mathcal{H}-  m_{2,5}\times \mathcal{J}\right)\notag\\
 &\quad +m_{1,5} \left(  m_{2,1} \times \mathcal{D}-  m_{2,2} \times \mathcal{G}+  m_{2,3} \times \mathcal{I}-  m_{2,4}\times \mathcal{J}\right).
\end{align}

Let $CF_{i,j}$ denote the cofactor of the entry in the $i$-th row and the $j$-th column of $A$. Then,
\begin{align}
CF_{1,3}&=\text{Det}\left(\begin{pmatrix}
A_{2}(q) & A_{1}(q)  & A_{6}(q) & A_{5}(q) & A_{4}(q) & A_{3}(q)\\
A_{3}(q) & A_{2}(q)  & A_{7}(q) & A_{6}(q) & A_{5}(q) & A_{4}(q)\\
A_{4}(q) & A_{3}(q)  & A_{1}(q) & A_{7}(q) & A_{6}(q) & A_{5}(q)\\
A_{5}(q) & A_{4}(q)  & A_{2}(q) & A_{1}(q) & A_{7}(q) & A_{6}(q)\\
A_{6}(q) & A_{5}(q)  & A_{3}(q) & A_{2}(q) & A_{1}(q) & A_{7}(q)\\
A_{7}(q) & A_{6}(q)  & A_{4}(q) & A_{3}(q) & A_{2}(q) & A_{1}(q)
\end{pmatrix}\right)=:\text{Det}\left(A_{\text{sub}}\right)\notag\\
\label{Minors A sub} &=A_2(q)\times U_{1,1}-A_1(q)\times U_{1,2}+A_6(q)\times U_{1,3}-A_5(q)\times U_{1,4}\notag\\
&\quad +A_4(q)\times U_{1,5}-A_3(q)\times U_{1,6},
\end{align}
where $U_{i,j}$ are the minors of the matrix $A_{\text{sub}}$.

We now employ the determinant function \eqref{Det function} to find $U_{1,j}$, $1\le j \le6$ in \eqref{Minors A sub}, which eventually give
\begin{align*}
 CF_{1,3}&= A_6^6(q)-5 A_2(q) A_3(q) A_6^4(q)-5 A_1(q) A_4(q) A_6^4(q)-5 A_5(q) A_7(q) A_6^4(q)\\
 &\quad
 +4 A_3(q) A_4^2(q) A_6^3(q)+4 A_1(q) A_5^2(q) A_6^3(q)+4 A_4(q) A_7^2(q) A_6^3(q)\\
 &\quad
 +4 A_1^2(q) A_2(q) A_6^3(q)+4 A_3^2(q) A_5(q) A_6^3(q)+8 A_2(q) A_4(q) A_5(q) A_6^3(q)\\
 &\quad
 +4 A_2^2(q) A_7(q) A_6^3(q)+8 A_1(q) A_3(q) A_7(q) A_6^3(q)-3 A_1(q) A_3^3(q) A_6^2(q)\\
 &\quad
 -3 A_2(q) A_5^3(q) A_6^2(q)-3 A_3(q) A_7^3(q) A_6^2(q)+6 A_2^2(q) A_3^2(q) A_6^2(q)\\
 &\quad
 +6 A_1^2(q) A_4^2(q) A_6^2(q)-9 A_3(q) A_4(q) A_5^2(q) A_6^2(q)+6 A_5^2(q) A_7^2(q) A_6^2(q)\\
 &\quad
 -9 A_1(q) A_2(q) A_7^2(q) A_6^2(q)-3 A_2^3(q) A_4(q) A_6^2(q)+3 A_1(q) A_2(q) A_3(q) A_4(q) A_6^2(q)\\
 &\quad
 -3 A_4^3(q) A_5(q) A_6^2(q)-9 A_1(q) A_2^2(q) A_5(q) A_6^2(q)-9 A_1^2(q) A_3(q) A_5(q) A_6^2(q)\\
 &\quad
 -3 A_1^3(q) A_7(q) A_6^2(q)-9 A_2(q) A_4^2(q) A_7(q) A_6^2(q)-9 A_3^2(q) A_4(q) A_7(q) A_6^2(q)\\
 &\quad
 +3 A_2(q) A_3(q) A_5(q) A_7(q) A_6^2(q)+3 A_1(q) A_4(q) A_5(q) A_7(q) A_6^2(q)+2 A_1(q) A_2^4(q) A_6(q)\\
 &\quad
 +2 A_3(q) A_5^4(q) A_6(q)+2 A_2(q) A_7^4(q) A_6(q)+4 A_2^2(q) A_4^3(q) A_6(q)\\
 &\quad
 -6 A_1(q) A_3(q) A_4^3(q) A_6(q)+4 A_4^2(q) A_5^3(q) A_6(q)+4 A_1^2(q) A_7^3(q) A_6(q)\\
 &\quad
 -6 A_4(q) A_5(q) A_7^3(q) A_6(q)+4 A_1^3(q) A_3^2(q) A_6(q)-2 A_2(q) A_3^2(q) A_4^2(q) A_6(q)\\
 &\quad
 +4 A_2^3(q) A_5^2(q) A_6(q)+10 A_1(q) A_2(q) A_3(q) A_5^2(q) A_6(q)-2 A_1^2(q) A_4(q) A_5^2(q) A_6(q)\\
 &\quad
 +4 A_3^3(q) A_7^2(q) A_6(q)-2 A_1(q) A_4^2(q) A_7^2(q) A_6(q)+10 A_2(q) A_3(q) A_4(q) A_7^2(q) A_6(q)\\
 &\quad
 -2 A_2^2(q) A_5(q) A_7^2(q) A_6(q)-4 A_1(q) A_3(q) A_5(q) A_7^2(q) A_6(q)\\
 &\quad
 -2 A_1^2(q) A_2^2(q) A_3(q) A_6(q)+2 A_3^4(q) A_4(q) A_6(q)-6 A_1^3(q) A_2(q) A_4(q) A_6(q)\\
 &\quad
 +2 A_1^4(q) A_5(q) A_6(q)-6 A_2(q) A_3^3(q) A_5(q) A_6(q)-4 A_1(q) A_2(q) A_4^2(q) A_5(q) A_6(q)\\
 &\quad
 +10 A_1(q) A_3^2(q) A_4(q) A_5(q) A_6(q)-4 A_2^2(q) A_3(q) A_4(q) A_5(q) A_6(q)\\
 &\quad
 +2 A_4^4(q) A_7(q) A_6(q)-6 A_1(q) A_5^3(q) A_7(q) A_6(q)-4 A_1(q) A_2(q) A_3^2(q) A_7(q) A_6(q)\\
 &\quad
 -2 A_3^2(q) A_5^2(q) A_7(q) A_6(q)-4 A_2(q) A_4(q) A_5^2(q) A_7(q) A_6(q)\\
 &\quad
 -6 A_2^3(q) A_3(q) A_7(q) A_6(q)+10 A_1(q) A_2^2(q) A_4(q) A_7(q) A_6(q)\\
 &\quad
 -4 A_1^2(q) A_3(q) A_4(q) A_7(q) A_6(q)+10 A_3(q) A_4^2(q) A_5(q) A_7(q) A_6(q)\\
 &\quad
 +10 A_1^2(q) A_2(q) A_5(q) A_7(q) A_6(q)-A_2(q) A_4^5(q)-A_4(q) A_5^5(q)-A_1(q) A_7^5(q)\\
 &\quad
 +2 A_1(q) A_2(q) A_3^4(q)+A_3^2(q) A_4^4(q)+A_1^2(q) A_5^4(q)+A_4^2(q) A_7^4(q)\\
 &\quad
 +2 A_3(q) A_5(q) A_7^4(q)-A_2^3(q) A_3^3(q)-A_1^3(q) A_4^3(q)-3 A_1(q) A_3^2(q) A_5^3(q)\\
 &\quad
 -3 A_2^2(q) A_3(q) A_5^3(q)+A_1(q) A_2(q) A_4(q) A_5^3(q)-A_5^3(q) A_7^3(q)-3 A_2(q) A_3^2(q) A_7^3(q)\\
 &\quad
 -3 A_2^2(q) A_4(q) A_7^3(q)+A_1(q) A_3(q) A_4(q) A_7^3(q)+A_1(q) A_2(q) A_5(q) A_7^3(q)\\
 &\quad
 +A_1^4(q) A_2^2(q)-3 A_1(q) A_2^3(q) A_4^2(q)+5 A_1^2(q) A_2(q) A_3(q) A_4^2(q)+A_3^4(q) A_5^2(q)\\
 &\quad
 -A_2^2(q) A_4^2(q) A_5^2(q)-2 A_1(q) A_3(q) A_4^2(q) A_5^2(q)-3 A_1^3(q) A_2(q) A_5^2(q)\\
 &\quad
 +5 A_2(q) A_3^2(q) A_4(q) A_5^2(q)+A_2^4(q) A_7^2(q)-3 A_3(q) A_4^3(q) A_7^2(q)\\
 &\quad
 -A_1^2(q) A_3^2(q) A_7^2(q)-2 A_2(q) A_3(q) A_5^2(q) A_7^2(q)+5 A_1(q) A_4(q) A_5^2(q) A_7^2(q)\\
 &\quad
 +5 A_1(q) A_2^2(q) A_3(q) A_7^2(q)-2 A_1^2(q) A_2(q) A_4(q) A_7^2(q)-3 A_1^3(q) A_5(q) A_7^2(q)\\
 &\quad
 +5 A_2(q) A_4^2(q) A_5(q) A_7^2(q)-2 A_3^2(q) A_4(q) A_5(q) A_7^2(q)-A_1^5(q) A_3(q)\\
 &\quad
 -3 A_1^2(q) A_3^3(q) A_4(q)-2 A_1(q) A_2^2(q) A_3^2(q) A_4(q)+2 A_2^4(q) A_3(q) A_4(q)\\
 &\quad
 -A_2^5(q) A_5(q)+2 A_1(q) A_4^4(q) A_5(q)+A_2(q) A_3(q) A_4^3(q) A_5(q)\\
 &\quad
 -2 A_1^2(q) A_2(q) A_3^2(q) A_5(q)-3 A_3^3(q) A_4^2(q) A_5(q)+A_1(q) A_2^3(q) A_3(q) A_5(q)\\
 &\quad
 +5 A_1^2(q) A_2^2(q) A_4(q) A_5(q)+A_1^3(q) A_3(q) A_4(q) A_5(q)-A_3^5(q) A_7(q)\\
 &\quad
 +2 A_2(q) A_5^4(q) A_7(q)-3 A_1^2(q) A_2^3(q) A_7(q)+A_1(q) A_2(q) A_4^3(q) A_7(q)\\
 &\quad
 +A_3(q) A_4(q) A_5^3(q) A_7(q)+5 A_1(q) A_3^2(q) A_4^2(q) A_7(q)-2 A_2^2(q) A_3(q) A_4^2(q) A_7(q)\\
 &\quad
 -3 A_4^3(q) A_5^2(q) A_7(q)-2 A_1(q) A_2^2(q) A_5^2(q) A_7(q)+5 A_1^2(q) A_3(q) A_5^2(q) A_7(q)\\
 &\quad
 +A_1^3(q) A_2(q) A_3(q) A_7(q)+2 A_1^4(q) A_4(q) A_7(q)+A_2(q) A_3^3(q) A_4(q) A_7(q)\\
 &\quad
 +A_1(q) A_3^3(q) A_5(q) A_7(q)+5 A_2^2(q) A_3^2(q) A_5(q) A_7(q)-2 A_1^2(q) A_4^2(q) A_5(q) A_7(q)\\
 &\quad
 +A_2^3(q) A_4(q) A_5(q) A_7(q)-15 A_1(q) A_2(q) A_3(q) A_4(q) A_5(q) A_7(q).
\end{align*}

Using \textit{Mathematica}, we expand the right side of the above identity in terms of $f_{1,7}\left(q^7\right)$, $f_{2,7}\left(q^7\right)$, and $f_{3,7}\left(q^7\right)$ using \eqref{A1 in f}--\eqref{A7 in f} and \eqref{R i j k}. The resulting expansion, with the help of the mathematica command \texttt{FullSimplify[Expression, Assumptions]}, under the equalities \eqref{Condition 1}--\eqref{Condition 2} can be simplified as
\begin{align}
\label{CF 13} CF_{1,3}&=7q^2\Big(-532544 f_{1,7}^{23}\left(q^7\right)+2822366 f_{1,7}^{22}\left(q^7\right) f_{2,7}^2\left(q^7\right)-9375040 f_{1,7}^{21}\left(q^7\right) f_{2,7}^4\left(q^7\right)\notag\\
&\quad
+21207130 f_{1,7}^{20}\left(q^7\right) f_{2,7}^6\left(q^7\right)-34041692 f_{1,7}^{19}\left(q^7\right)f_{2,7}^8\left(q^7\right)+39647716 f_{1,7}^{18}\left(q^7\right) f_{2,7}^{10}\left(q^7\right)\notag\\
&\quad
-34010032 f_{1,7}^{17}\left(q^7\right) f_{2,7}^{12}\left(q^7\right)+21762764 f_{1,7}^{16}\left(q^7\right) f_{2,7}^{14}\left(q^7\right)-10688908 f_{1,7}^{15}\left(q^7\right) f_{2,7}^{16}\left(q^7\right)\notag\\
&\quad
+4393575 f_{1,7}^{14}\left(q^7\right) f_{2,7}^{18}\left(q^7\right)-1624042 f_{1,7}^{13}\left(q^7\right) f_{2,7}^{20}\left(q^7\right)+347825 f_{1,7}^{12}\left(q^7\right) f_{2,7}^{22}\left(q^7\right)\notag\\
&\quad
+133384 f_{1,7}^{11}\left(q^7\right) f_{2,7}^{24}\left(q^7\right)-115269 f_{1,7}^{10}\left(q^7\right) f_{2,7}^{26}\left(q^7\right)+17154 f_{1,7}^9\left(q^7\right) f_{2,7}^{28}\left(q^7\right)\notag\\
&\quad
+3047 f_{1,7}^8\left(q^7\right) f_{2,7}^{30}\left(q^7\right)+36 f_{1,7}^7\left(q^7\right) f_{2,7}^{32}\left(q^7\right)+54691 f_{1,7}^{26}\left(q^7\right) f_{3,7}^2\left(q^7\right)\notag\\
&\quad
-2174 f_{1,7}^{29}\left(q^7\right) f_{3,7}^4\left(q^7\right)+15 f_{1,7}^{32}\left(q^7\right) f_{3,7}^6\left(q^7\right)\Big).
\end{align}

Thus, solving \eqref{Matrix Eq 1}, we have
\begin{align}
\label{Witness 49n+19 Eve 4} N_{2}(q)&= \dfrac{ D^4(q) }{f_{49}^4}\times \dfrac{CF_{1,3}}{\text{Det}(A)}.
\end{align}

Now, to give $N_2(q)$ its final form, we find Det$(A)$. For that, first $V:=\left(v_{i,j}\right)_{1\le i,j\le 7}$ be a square matrix of order 7. Then, for integers $\ell\ge1$, it can be shown by induction on $\ell$ that the entry in the $i$-th row and $j$-th column of the matrix product $V^{\ell}$ is given by
\begin{align}
\label{kth power entry formula}V^\ell=\left(\sum_{t_1}^7\sum_{t_2}^7\sum_{t_3}^7\cdots\sum_{t_{\ell-1}}^7v_{i,t_{\ell-1}}v_{t_{\ell-1},t_{\ell-2}}v_{t_{\ell-2},t_{\ell-3}}\cdots v_{t_{2},t_{1}}v_{t_{1},t_{j}}\right)_{1\le i,j\le 7}. 
\end{align}

Now, for the matrix
\begin{align*}
W:=\begin{pmatrix}
f_{1,7}\left(q^7\right) & 0 & q^5f_{3,7}\left(q^7\right) & 0 & 0 & -q^2 & -q f_{2,7}\left(q^7\right)\\
-q f_{2,7}\left(q^7\right) & f_{1,7}\left(q^7\right) & 0 & q^5f_{3,7}\left(q^7\right) & 0 & 0 & -q^2\\
-q^2 & -q f_{2,7}\left(q^7\right) & f_{1,7}\left(q^7\right) & 0 & q^5f_{3,7}\left(q^7\right) & 0 & 0\\
0 & -q^2 & -q f_{2,7}\left(q^7\right) & f_{1,7}\left(q^7\right) & 0 & q^5f_{3,7}\left(q^7\right) & 0\\
0 & 0 & -q^2 & -q f_{2,7}\left(q^7\right) & f_{1,7}\left(q^7\right) & 0 & q^5f_{3,7}\left(q^7\right)\\
q^5f_{3,7}\left(q^7\right) & 0 & 0 & -q^2 & -q f_{2,7}\left(q^7\right) & f_{1,7}\left(q^7\right) & 0\\
0 & q^5f_{3,7}\left(q^7\right) & 0 & 0 &-q^2 & -q f_{2,7}\left(q^7\right)& f_{1,7}\left(q^7\right)
\end{pmatrix},
\end{align*}
using \eqref{kth power entry formula}, we obtain
\begin{align}
\label{Matrix Equality} W^4=A.
\end{align}
Also by computing the determinant of $W$, we find that
\begin{align*}
\text{Det}(W)&=f_{1,7}^7\left(q^7\right)-q^7 \Big(f_{2,7}^7\left(q^7\right)+7 f_{1,7}\left(q^7\right) f_{2,7}^5\left(q^7\right)+14 f_{1,7}^2\left(q^7\right) f_{2,7}^3\left(q^7\right)\\
&\quad -7 f_{1,7}^4\left(q^7\right) f_{3,7}\left(q^7\right) f_{2,7}^2\left(q^7\right)+7 f_{1,7}^3\left(q^7\right) f_{2,7}\left(q^7\right)-7 f_{1,7}^5\left(q^7\right) f_{3,7}\left(q^7\right)\Big)\\
&\quad-q^{14} \Big(-7 f_{1,7}\left(q^7\right) f_{3,7}^2\left(q^7\right) f_{2,7}^4\left(q^7\right)-7 f_{1,7}^2\left(q^7\right) f_{3,7}^2\left(q^7\right) f_{2,7}^2\left(q^7\right)\\
&\quad +14 f_{1,7}\left(q^7\right) f_{3,7}\left(q^7\right) f_{2,7}\left(q^7\right)-14 f_{1,7}^3\left(q^7\right) f_{3,7}^2\left(q^7\right)+7 f_{3,7}\left(q^7\right) f_{2,7}^3\left(q^7\right)+1\Big)\\
&\quad -q^{21} \Big(7 f_{1,7}^2\left(q^7\right) f_{2,7}\left(q^7\right) f_{3,7}^4\left(q^7\right)-14 f_{2,7}^2\left(q^7\right) f_{3,7}^3\left(q^7\right)-7 f_{1,7}\left(q^7\right) f_{3,7}^3\left(q^7\right)\Big)\\
&\quad -7 q^{28} f_{2,7}\left(q^7\right) f_{3,7}^5\left(q^7\right)+q^{35} f_{3,7}^7\left(q^7\right)\\
&=\prod_{j=0}^6\left(f_{1,7}\left(q^7\right)- \zeta ^jq f_{2,7}\left(q^7\right)-\zeta ^{2 j}q^2+\zeta ^{5 j}q^5 f_{3,7}\left(q^7\right)\right),
\end{align*}
which, by \eqref{7 Dissection f1}, gives
\begin{align}
\label{Witness 49n+19 Eve 5} \text{Det}(W)&=\dfrac{\prod_{j=0}^6\left(\zeta^j q;\zeta^j q\right)_\infty}{f_{49}^7}=\dfrac{D(q)}{f_{49}^7}.
\end{align}
Therefore, employing \eqref{Witness 49n+19 Eve 4}, \eqref{Matrix Equality}, and \eqref{Witness 49n+19 Eve 5}, we have
\begin{align}
\label{Witness 49n+19 Eve 6} N_2(q)&=f_{49}^{24}\times CF_{1,3}.
\end{align}

Now, from \eqref{Witness 49n+19 Eve 1.5}, extracting the terms that contain $q^{7n+2}$, and then using \eqref{Witness 49n+19 Eve 3} and \eqref{Witness 49n+19 Eve 6}, we obtain
\begin{align*}
\sum_{n=0}^{\infty}p_{1}(7n+2)q^{7n+2}&=\dfrac{f_7^3\times N_2(q)}{D^4(q)}=\dfrac{ f_{49}^{28}}{f_7^{29}}\times CF_{1,3}.
\end{align*}
The above identity, with the help of \eqref{CF 13}, evidently gives
\begin{align}
\label{Final 1} \sum_{n=0}^{\infty}p_{1}(7n+2)q^{n}&=7\dfrac{f_{7}^{28}}{f_{1}^{29}}\Bigg(\sum_{j=0}^{16}\alpha_{19,j} f_{1,7}^{23-j}(q)f_{2,7}^{2j}(q)+\sum_{j=0}^{2}\beta_{19,j} f_{1,7}^{3j+26}(q)f_{3,7}^{2j+2}(q)\Bigg),
\end{align}
where $\alpha_{19,j}$ and $\beta_{19,j}$ are given in Table \ref{Table 2}. One can similarly find that
\begin{align}
\label{Final 2} \sum_{n=0}^{\infty}p_{2}(7n+2)q^{n}&=\dfrac{f_{7}^{56}}{f_{1}^{57}}\Bigg(\sum_{j=0}^{33}\gamma_{19,j} f_{1,7}^{47-j}(q)f_{2,7}^{2j+1}(q)+\sum_{j=0}^{6}\delta_{19,j} f_{1,7}^{3j+49}(q)f_{3,7}^{2j+1}(q)
\Bigg),
\end{align}
where $\gamma_{19,j}$ and $\delta_{19,j}$ are given in Tables \ref{Table 3} and \ref{Table 4}. Thus, \eqref{Witness 49n+19 Eve 1}, \eqref{Final 1}, and \eqref{Final 2} together imply \eqref{Witness 49n+19}.\qed

\section{Concluding remarks}
The method of this paper is an extension of \cite{GJS22}. It seems that this extension can be applied for finding generating functions for the 2-color partition functions as well. Let $p_{1,r}(n)$ denote the number of partitions, where the parts that are multiple of $r$ can appear in two colors and the others appear in one color. Then, it would be interesting to build further on this extension, in quest of finding the generating functions of the following congruences. For all $n\ge 0$, we have
\begin{align*}
p_{1,7}(25n+17)&\equiv 0\pmod{5},\\
p_{1,17}(25n+7)&\equiv 0\pmod{5},\\
p_{1,2}(49n+t)&\equiv 0\pmod{7},\quad t\in\{15,29,36,43\},\\
p_{1,4}(49n+t)&\equiv 0\pmod{7},\quad t\in\{11,25,32,39\}.
\end{align*}

\section*{Acknowledgements}
The third author was partially supported by the institutional fellowship for doctoral research from Tezpur University, Napaam, India. The author thanks the funding institution.

\section{Appendix}
The values of $\alpha_{r,j}$,  $\beta_{r,j}$, $\gamma_{r,j}$, and $\delta_{r,j}$ for $r\in\{19,33,40\}$ are given in the following tables.

\begin{longtable}{c c c c c c c} 
\caption{Values of $\alpha_{19,j}$, $\alpha_{33,j}$, $\alpha_{40,j}$, $\beta_{19,j} $, $\beta_{33,j} $, and $\beta_{40,j} $}
\label{Table 2}
\\
 \hline
 $j$  & $\alpha_{19,j}$ & $\alpha_{33,j}$ & $\alpha_{40,j}$ & $\beta_{19,j} $ & $\beta_{33,j} $ & $\beta_{40,j} $  \T\B\\  \hline
0 & -532544  & 677424 & 1726428 & 54691 & -78156 & -288131 \T\B\\
1 & 2822366 & -3247912 & -6229360  & -2174 & 3731 & 24263\T\B\\
2 & -9375040  & 9923453 & 14952497 & 15 & -36 & -662\T\B\\
3 & 21207130 & -20610156 & -24984853 & & & 2 \T\B\\
4 & -34041692 & 30360672 & 29917296 \T\B\\
5 & 39647716  & -32527768 & -26244924 \T\B\\
6 & -34010032 & 25602765 & 17026404 \T\B\\
7 & 21762764 & -14863836 & -7970560    \T\B\\
8 & -10688908 & 6287362 & 2294400   \T\B\\
9 & 4393575 & -1603020 & -25088    \T\B\\
10 & -1624042 & -208775 & -232970    \T\B\\
11 & 347825 & 412284 & -23794    \T\B\\
12 & 133384 & -116340 & 79116    \T\B\\
13 & -115269 & -24992 & -18124    \T\B\\
14 & 17154 & 12689 & -1889    \T\B\\
15 & 3047 & 624 & -15    \T\B\\
16 & 36 & 2 &     \T\B\\
\hline
\end{longtable}

\vspace*{.4cm}
\begin{longtable}{c c c c} 
\caption{Values of $\gamma_{19,j}$, $\gamma_{33,j}$, and $\gamma_{40,j}$}
\label{Table 3}
\\
 \hline
 $j$  & $\gamma_{19,j}$ & $\gamma_{33,j}$ & $\gamma_{40,j}$   \T\B\\  \hline
0 & 402306358809680  & -511808960773941 & 217368522907096  \T\B\\
1 & -2452400077267696 & 2965591662257740 & -1414744004689296  \T\B\\
2 & 11508628880956008  & -13276581304101572 & 7003958111796106 \T\B\\
3 & -42812168519266872 & 47229868718929840 & -27236672553224208  \T\B\\
4 & 129010608470924891 & -136309966462773336 & 85183453951673440 \T\B\\
5 & -320094377783048275  & 324206604364212679 & -218078089624300176 \T\B\\
6 & 662096611024442804 & -643106461893717549 & 463167416613039502 \T\B\\
7 & -1152641373328950740 & 1073688904333614081 & -824470260114451208    \T\B\\
8 & -1701256373500039024 & -1519299185881589891 & 1239707905766475508   \T\B\\
9 & -2140744449656720288 & 1831807622321672560 & -1583999740415812448    \T\B\\
10 & 2306190155043227112 & -1889319680591388592 & 1727526156193273852    \T\B\\
11 & -2133495352783794352 & 1671734309523742924 & -1613451148734963728    \T\B\\
12 & 1698596696577538715 & -1271563604540253748 & 1293486609171772752    \T\B\\
13 & -1165486377041496475 & 832518742908855716 & -891509417208035472    \T\B\\
14 & 689787793489082244 & -469527310455304100 & 528767666556150340    \T\B\\
15 & -352299212195429308 & 228140690684456640 & -269998559740321136    \T\B\\
16 & 155261889913586318 & -95458813361311800 &  118694926805304860   \T\B\\
17 & -58943833568644310 & 34418612386331232 & -44966370616410672  \T\B\\
18 & 19182508457950588 & -10828615634044104 & 14774887071736473   \T\B\\
19 & -5384836913311692 & 3128357985329824 & -4261261546744684   \T\B\\
20 & 1485787932113975 & -872227173968528 &  1006585751346434  \T\B\\
21 & -595456845326039 & 153762390276994 &  -46403764763292  \T\B\\
22 & 315164225806868 & 85318408921842 & -164330600769945  \T\B\\
23 & -113573559663364 & -98502880754070 & 109458165551552   \T\B\\
24 & 1145368735660 & 39501282070410 &  -23392382620184   \T\B\\
25 & 19646580077884 & -2980862040844 &  -8600187380856   \T\B\\
26 & -7875251361048 & -3094740306164 &  6171021799749   \T\B\\
27 & 605267542600 & 888128348156 &  -933784846260   \T\B\\
28 & 284714584631 & 19467043676 &   -154321594618  \T\B\\
29 & -36921650599 & -21811591791 &  34529888556   \T\B\\
30 & -5358304792 & -1361107667 &   3522219017  \T\B\\
31 & -146117528 & -20367821 &  78187824   \T\B\\
32 & -935474 & -61393 &  394676   \T\B\\
33 & -726 & -8 &  192   \T\B\\
\hline
\end{longtable}

\begin{longtable}{c c c c} 
\caption{Values of $\delta_{19,j}$, $\delta_{33,j}$, and $\delta_{40,j}$}
\label{Table 4}
\\
 \hline
 $j$  & $\delta_{19,j}$ & $\delta_{33,j}$ & $\delta_{40,j}$   \T\B\\  \hline
0 & -48785839561732  & 65651447803543 & -24286463387760  \T\B\\
1 & -4126733575196 & -5923164680395 & 1849125804685  \T\B\\
2 & -223442523301  & 346515196113 & -87052955380 \T\B\\
3 & 6787708205 & -11611129970 & 2180656118  \T\B\\
4 & -93729948 & 183167294 & -22768388 \T\B\\
5 & 402356  & -964514 & 61713 \T\B\\
6 & -192 & 726 & -8 \T\B\\
\hline
\end{longtable}


\end{document}